\documentclass[12pt]{amsart}
\usepackage{times}
\usepackage{amscd}
\usepackage[latin1]{inputenc}
\usepackage{latexsym}
\usepackage{amsmath,amssymb,amsthm,amsfonts}
\usepackage[dvips]{graphicx}
\newtheorem{teor}{Theorem}[section]
\newtheorem{lemma}[teor]{Lemma}
\newtheorem{propos}[teor]{Proposition}

\theoremstyle{definition}
\newtheorem{defin}[teor]{Definition}
\newtheorem{osserv}[teor]{Remark}
\newtheorem{es}[teor]{Example}

\newcommand{\C}{\mbox{ $ \mathbb{C} $ }}

\newcommand{\Diff}{\mbox{\textrm{Diff}}}

\theoremstyle{remark}  
\begin{document}
\title[Upper-bound for the number of RP  curves of a map tangent to identity ]
{Upper-bound
for the number of robust parabolic curves for a class of maps
tangent to identity }
\author[Francesco Degli Innocenti  \and Chiara Frosini]
{Francesco Degli Innocenti $^\spadesuit$ \and Chiara
Frosini$^\clubsuit$}

\address{$\begin{aligned}^\spadesuit & \hbox{Dipartimento di Matematica, Universit\`a di
Pisa, Largo Bruno Pontecorvo 5, 56127, Pisa.}\\ ^\clubsuit &
\hbox{Dipartimento di Matematica, Universit\`a di Firenze, viale
Morgagni 67/A, 50134, Firenze.}\\
 & \hspace*{3cm}\hbox{\textbf{\emph{e-mail address:}} \textit{degliinno@dm.unipi.it \and frosini@math.unifi.it}}\end{aligned}$}

\maketitle
\section{Introduction}
The Leau-Fatou flower theorem \cite{CG} completely describes the dynamic behavior of $1-$di\-men\-sio\-nal maps tangent to the identity. In
dimension two Hakim \cite{H} and Abate \cite{AB} proved that if $f$ is a holomorphic map tangent to the identity in $\mathbb{C}^2$ and $\nu(f)$
is the degree of the first non vanishing jet of $f-Id$ then there exist  $\nu(f)-1$  robust parabolic curves (RP curves for short), namely
attractive petals at the origin which survive under by blow-up (see \cite{AT} and Section \ref{Robust parabolic curves}) .
The set of the exponential of holomorphic vector fields (of order greater than or equal to two), $\Phi_{\geq 2}(\mathbb{C}^2,0)$, is
dense in the space of germs of maps tangent to the identity.

In this paper we give an upper-bound for the number of robust parabolic curves of $f\in \Phi_{\geq 2}(\mathbb{C}^2,0) .$

\begin{teor}\label{teorema:nostro}
Let $f=(f_{1}, f_{2})\in \Phi_{\geq 2}(\mathbb{C}^2,0)$ be a non-dicritical holomorphic map.
Set $\eta(f):=\max\{ord(f_{1}-Id),
ord(f_{2}-Id)\}$ and $\mu(f)$ the Milnor number of $f.$ Then the number of RP curves is at most $$(\mu(f)+1)(\eta^2(f)-\eta(f)).$$
\end{teor}
In the dicritical  case in \cite{BR} Bracci proved that $f$ is dicritical if and only if it
 has infinitely many parabolic curves. Here we show that the
parabolic curves at a dicritical point are indeed robust ones:
\begin{propos}
Let $f=(f_{1}, f_{2})\in \Phi_{\geq 2}(\mathbb{C}^2,0)$ be a dicritical holomorphic map. Then there exist infinitely many  RP curves.
\end{propos}

A new approach to study dynamics of germs of mappings tangent to the identity has been introduced in \cite{BR}, \cite{ABT} by Abate, Bracci and Tovena.
Their idea is to study discrete dynamics using families of vector fields whose flows are approximations at the first order of $f$.
This technique turns out to be is very useful.

In case $f$ is the exponential of an holomorphic vector field we
notice a very strict relation between the dynamics of the map and
the dynamics of the vector field. So, if the diffeomorphism $f,$
is such that there exists a vector field $X$ such that $
\exp(X)=f,$ then it turns out that the RP curves are
"geometrically" determined by the fact they lay in a separatrix of
$X.$ Thus we can use a result of Corral and Fernandez Sanchez
\cite{CFS} concerning the upper-bound of the number of
separatrices of $X$ to estimate the number of RP curves  .

\section{Exponential of a vector field}
Let $f\in \Phi_{\geq 2}(\mathbb{C}^2,0)$ be a holomorphic map
tangent to the identity in $\mathbb{C}^2$ then there exists a
vector field $X$ such that $ \exp(X)=F.$ The two objects, $f$ and
$X$, are related by the following classical result:

\begin{propos}\label{propos:formalflow}
Let $X$ be a germ of holomorphic vector in $(\mathbb{C}^2,0).$ Then its time one flow can be written has:
\begin{equation}\label{eq:conv}
f^t(x,y)=\left( x+\sum_{n=1}^\infty \frac{t^n}{n!}X^n.x, y+ \sum_{n=1}^{\infty}\frac{t^n}{n!}X^n.y \right ),
\end{equation}
where $X^n.x$ is defined by $X$ applied to $X^{n-1}.x$ and $X.x$ is the application of $X$ to $x$.
\end{propos}

\begin{osserv}
If $X$ is an holomorphic vector field in $(\mathbb{C}^2,0)$ then its time one map is given by:
\begin{equation}\label{eq:formal}
\exp(\hat{X})=\left(x+\sum_{n=1}^{\infty}\frac{1}{n!}\hat{X}^n.x,
y+\sum_{n=1}^{\infty}\frac{1}{n!}\hat{X}^n.y \right).
\end{equation}
\end{osserv}

\begin{osserv}
If $f$ is a holomorphic map tangent to the identity in $\mathbb{C}^2$ at the origin
then the associated vector field has order of singularity at $0$ greater or equal to $2.$
\end{osserv}

\section{Robust parabolic curves}\label{Robust parabolic curves}
In this section we give the definition
of parabolic and robust parabolic curves for a germ of holomorphic map tangent to the identity at the origin in $\mathbb{C}^2$
and we study the relationship between the separatrices of the vector field associated to $f$ and these curves.

\begin{defin}
A \textbf{parabolic curve} for $f\in \Diff(\mathbb{C}^2,0)$ at the origin is an injective holomorphic map $ \varphi: \Delta \rightarrow
\mathbb{C}^2$ satisfying the following properties:
\begin{enumerate}
\item $\Delta$ is a simply connected domain in $\mathbb{C}$ with $0\in \partial \Delta;$ \item $\varphi$ is continuous at the origin, and
$\varphi(0)=0;$ \item $\varphi(\Delta)$ is invariant under $f,$ and $(f\mid_{\varphi(\Delta)})^n \rightarrow 0$ as $n\rightarrow \infty.$
\end{enumerate}
\end{defin}

The idea of robust parabolic curve was first introduced by Abate and Tovena in \cite{AT}:
\begin{defin}\label{def:curve paraboliche robuste}
A \textbf{robust parabolic curve} is a parabolic curve that satisfies the following conditions:
\begin{enumerate}
\item we can blow-up $\varphi$ at level $h$ for any $h\geq 1,$ \item there is a formal power series $Q\in (\mathbb{C}[[x]])^2$ such that for
every $h\geq 1$ there is $r_h>0$ such that $ \varphi-Q_h = O(\zeta^{h+1})$ in $ \Delta_{r_h},$ where $Q_h$ denotes the truncation at degree $h$
of $Q.$
\end{enumerate}
\end{defin}
\begin{osserv}
Essentially when we say condition $(1)$ we mean that the strict transform of the parabolic curve is also a parabolic curve, but for a clear
definition of "blow-up at level $h$" we refer to \cite{AT}.
\end{osserv}

The geometric meaning of Definition \ref{def:curve paraboliche robuste} is clarified by the following
proposition:
\begin{propos}\label{propos:sepfor}
Let $f\in \Phi_{\geq 2}(\mathbb{C}^2,0)$ be a holomorphic map and let $X$ be vector field such that $ \exp(X)=f$.
Let $\varphi$ be a robust parabolic curve. Then $\varphi$ is
contained in a
separatrix of $X.$ Conversely in every
separatrix of $X$ there exists at least one RP curve for $f$.
\end{propos}
\begin{proof}
Let be $p\in \varphi(\Delta).$
Since $$
\exp(X)=f.
$$
then the orbit $\{f^n(p)\}$ is contained in a separatrix, $S,$ of $X.$
We have only to prove that every orbit generated by a generic point  $q\in\varphi(\Delta)$ stays inside $S$.
By contradiction we
can find two orbits that converge to zero living in two different separatrices , say $S_{1}$ and  $S_{2}$.
Let $l_{1}(x,y)$ and $l_{2}(x,y)$ be (respectively) the local expressions of $S_{1}$ and  $S_{2}$
and let $h$ be the order of the first non zero jet of $l_{1}-l_{2}.$
If we blow-up the vector field $h$ times then, by property $(1)$ of the definition of RP curves,
we have that the two orbits converge to zero
with two different directions and this contradicts property $(2)$ of Definition \ref{def:curve paraboliche robuste}.

Let prove the converse.
Let be $S$ a  separatrix and let $ y=\varphi(x)=x^{\frac{p}{q}}+ \cdots$ be its expression in
Puiseux series.
\begin{osserv}\label{osserv:exponent}
We can suppose $p$ and $q$ prime each other and that $\frac{p}{q}\geq 1,$ if this is not true we can choose the parametrization of the
separatrix in the form $x=\psi(y)$ which satisfies the required condition.
\end{osserv}

Let now make the following change of variables:
$$
\left \{
\begin{aligned}
u=&x\\
v=&y-\varphi(x)
\end{aligned}
\right.
$$
The vector field in the new coordinates is:
$$
\left\{
\begin{aligned}
\dot{u}=& A(u,v+\varphi(u))\\
\dot{v}=&B(u,v+\varphi(u))-\dot{\varphi}(u)A(u,v+\varphi(u))
\end{aligned}
\right.
$$
If we compute the exponential of this new vector field restricted to the separatrix $\{v=0\}$we find:
$$
\exp(A(u,\varphi(u))\frac{\partial}{\partial u})
$$
Let us make the change of variables:
$$
u=z^q
$$
and then the first component of the vector field is:
\begin{equation}\label{eq:newfield}
\dot{z}=\frac{A(z^q,\varphi(z^q))}{qz^{q-1}}.
\end{equation}
By Remark \ref{osserv:exponent} the left-hand side of  \eqref{eq:newfield} is expressed as a power series.
 Now if we take the
exponential of this new vector field we find a map tangent to the identity conjugated to the original one given by
$ (z,w)\mapsto (z+ z^h + \cdots,w).$ By the Leau-Fatou
Theorem (\cite{CG}) we get the assertion.
\end{proof}
As a consequence of this last result we easily prove
the existence of RP curves for map in $\Phi_{\geq 2}(\mathbb{C}^2,0)$ \cite{AB},\cite{H}, \cite{AT}.
\begin{propos}
Let $f\in \Phi_{\geq 2}(\mathbb{C}^2,0)$ be a holomorphic map tangent to the identity in $\C^2.$
Then there exists at least one RP curve.
\end{propos}

\section{Non-dicritical case}
Proposition \ref{propos:sepfor} shows that the RP curves live inside separatrices of the associated vector
field. So
the idea  is to estimate the number of separatrices of the vector field and then  the number of RP curves
inside a separatrix.
In \cite{CFS} Corral and Fernandez Sanchez find the optimal estimation of the number of separatrices by the Milnor number of $X$ \cite{BPVV}.
\begin{propos}\cite{CFS}
Let  $X$ be an holomorphic vector field in $\mathbb{C}^2,$ singular at the origin.
Let $S$ be  the curve determined by all the separatrices trough the
origin. Let  $r_0 (S)$ be  the number of the irreducible components of $S.$ Then:
\begin{equation}
r_0(S)\leq \mu_0(X)+1,
\end{equation}
where $\mu_0$ is the Milnor number of $X$ at the origin.
\end{propos}
The proof of this proposition can be found in \cite{CFS}.
In order to express the previous estimation in terms of invariants of $f$ we introduce the intersection multiplicity \cite{AB}.
\begin{propos}
Let $f\in \Phi_{\geq 2}(\mathbb{C}^2,0)$ be a holomorphic map tangent to the identity in $\C^2$ and let $X$ be the associated vector field.
Then the Milnor number of $X,$ at the origin, is equal to the
intersection multiplicity of $f-Id.$
\end{propos}
\begin{proof}
Let observe that, if
$$X=A(x,y)\frac{\partial}{\partial x}+B(x,y)\frac{\partial}{\partial y},$$
then the Milnor number of $X$ is equal to the intersection multiplicity at the origin \cite{BPVV}.
Now observe that the intersection
multiplicity of two function $g(x,y)$ and $h(x,y)$ depends only on the first non zero jet
of $g$ and on the lowest exponent of the Puiseux
parametrization of $h$ \cite{Gr}. The Newton-Puiseux polygon shows that the lowest exponent of the parametrization
 depends only on the part of the polygon
determined by the first non zero jet of the function \cite{Cas}.
This concludes the proof because the lowest non zero jets of $f-Id,$
 $A$ and $B$ are the same.
\end{proof}

 Now we can start with the second part of the proof i.e. the estimation of  the number of RP curves
 that are  in a  separatrix.
 Proceeding as in Proposition \ref{propos:sepfor}, we can conjugate the restriction of $f$ to the separatrix
 to a map of the kind $(z,w)\mapsto(z+z^h+\cdots, w).$

 We have now to estimate the exponent $h.$ An easy computation
shows that the exponent $h$ is the lowest degree of the expression of
$\frac{A(z^q,\varphi(z^q))}{qz^{q-1}}.$
The same computation proves that the order of $A(x,y)$ is the same as the order, $\nu_{1},$
 of $f_1 - Id.$
Then
$$
\frac{A_{\nu_1}(z^q,\varphi(z^q))}{qz^{q-1}}=\sum_{i+j=\nu_1} z^{qi}\varphi(z^q)^j z^{1-q},
$$
so the lowest degree is:
\begin{equation}\label{eq:parab}
qi+pj+1-j,
\end{equation}
for $0\leq i,j\leq \nu_1.$
Let us maximize the quantity \eqref{eq:parab}.
 According to the cases $ p,q>0$ and $p,q<0$ and by the assumption
$\frac{p}{q}\geq 1$ we have that:
$$
qi+pj+1-j\leq \nu_1 p+ 1 - q\leq \nu_1 p,
$$
where the last inequality holds because $q\geq 1.$
By Remark \ref{osserv:exponent} the number of RP curves in $S$ is estimated from above by:
$$
\max\{\nu_1, \nu_2\}p.
$$
This estimation depends on $p$ and $q$ and then on the  particular separatrix.
It is possible to improve this result removing the dependence on the separatrix in the following way.
Since:
$$
\frac{dy}{dx}=\frac{B(x,y)}{A(x,y)}
$$
and replacing $y=\varphi(x)=x^k+\cdots,$ we find that:
$$
k=\frac{i_1-i_1+1}{j_2-j_1+1},
$$
with $i_1+j_1=\nu_1$ and $i_2+j_2=\nu_2.$
So we easily find:
$$
k\leq \frac{\nu_2 -1}{\nu_1 -1}.
$$
Then $p\leq (\nu_2 -1)$ and $q\leq (\nu_1 -1)$ because $p$ and $q$ are prime each other.

This proves the main Theorem \ref{teorema:nostro}.

\section{Dicritical case}
In this section we estimate the number of
RP curves in the dicritical case. We briefly recall
 the definition of dicritical singularity for maps and for
vector fields.
\begin{defin}
Let $X$ be an holomorphic  vector field and let $p$ be a singularity.
We say that $p$ is a dicritical singularity if the exceptional divisor is not
invariant for the strict transform of $X.$
\end{defin}
For the case of maps tangent to the identity in $\mathbb{C}^2$ we have the following definitions:
\begin{defin}\cite{BR}
Let $S$ be an irreducible curve in $\C^2$ and $f:\C^2 \rightarrow \C^2$ be an holomorphic map such that
$f\mid_S=Id_S$ and $ f \neq Id.$ Let $\mathcal{I}(S)_p \subset \mathcal{O}_p$ be the
ideal of germs vanishing on $S.$ We say that $f$ is
tangential on $S$ at $p$ if for a defining function $l$ of $S$ at $p$:
\begin{equation*}
\frac{l\circ f - l}{l^T} \equiv 0 \hbox{\hspace{0.5cm}  mod $\mathcal{I}(S)_p$,}
\end{equation*}
where $T$ is the order of $f$ on $S$ at $p.$
\end{defin}
\begin{defin}\cite{BR}\cite{AB}\cite{BT}
Let $f$ be a holomorphic map tangent to the identity in $\mathbb{C}^2.$ We say that $0$ is dicritical for $f$ if the
blow-up of $f$ is non-tangential on the exceptional divisor.
\end{defin}
By the theorem on existence and uniqueness of solutions of ordinary differential equations \cite{Se} we have:
\begin{propos}\label{prop:inf separ}
Let $X$ be an holomorphic vector field and $p$ a dicritical singularity, then there exist infinitely many  separatrices trough $p$.
\end{propos}
In the discrete case the following result holds:
\begin{propos}\cite{AB}\cite{BR}
Let $f$ be a map tangent to the identity in $\C^2$. If $0$ is dicritical then there exists infinitely many parabolic curves for
$f.$
\end{propos}
The last two propositions suggest that there exist a strict relationship
between the two notions of dicriticity:
\begin{propos}\label{prop:dicritico}
Let  $f\in \Phi_{\geq 2}(\mathbb{C}^2,0)$ be a map tangent to the identity in $\C^2$ and $X$ the vector field
such that $exp(X)=f.$ Then $f$ is dicritical in
$0$ if and only if $X$ is dicritical in $0.$
\end{propos}
Before  proving the previous proposition
we need the following lemma, that is a direct application of Proposition $4.2.4$ (pg. 267) of
\cite{AMR}:
\begin{lemma}\label{lemma: marsed}
Let  $f\in \Phi_{\geq 2}(\mathbb{C}^2,0)$ be a map tangent to the identity in $\C^2$ and set $\tilde{f}$  the blow-up of $f.$
Let  $X$ be the vector field
associated to $f.$ We have that:
$$
\tilde{f}=exp(\pi^\ast X),
$$
where $\pi^\ast X$ is the pull-back of $X$ by the blow-up map $\pi.$
\end{lemma}
\begin{osserv}
This lemma says that the relationship between maps and vector fields is
preserved under blow-up if, instead of saturating the vector field, we divide it only
by the first power of a local expression of the exceptional divisor.
\end{osserv}
\begin{es}
This example shows that the equality $\tilde{f}= \exp(\tilde{X})$ is not generally true. Let:
$$
f(x,y)=(x+\sum_{k\geq 1} x^{k+1}, y+\sum_{k\geq 1}y^{k+1}).
$$
The associated vector field $X$ is:
$$
X=x^2\frac{\partial}{\partial x}+y^2\frac{\partial}{\partial y},
$$
i.e.
$$
\exp(X)=f.
$$
The blow-up of $f$ and $X$, in the chart for which $\pi(u,v)=(u,uv),$ are:
$$
\begin{aligned}
\tilde{f}&=(u+\sum_{k\geq 1} u^{k+1}, \frac{v+\sum_{k\geq 1}u^k v^{k+1}}{1+\sum_{k\geq 1} u^k});\\
\tilde{X}&=u\frac{\partial}{\partial u}+ (v^2-v)\frac{\partial}{\partial v}.
\end{aligned}
$$
So we easily prove that:
$$
\exp(\tilde{X})=(eu, \cdots)\neq \tilde{f}.
$$
But if we take:
$$
\bar{X}=\pi^\ast(X)=u^2\frac{\partial}{\partial u}+ u(v^2-v)\frac{\partial}{\partial v},
$$
we see that $\exp(\bar{X})=\tilde{f}.$ We can observe that generally $\exp(\tilde{X})$
is not tangent to the identity because the vector field has linear part.
\end{es}
\begin{proof}
If $X$ is dicritical then by Proposition \ref{prop:inf separ} there exists infinitely many separatrices and then, by Proposition \ref{propos:sepfor},
$f$ admits infinitely many RP curves. So $f$, by Theorem \ref{teorema:nostro}, has to be dicritical.
Let prove the converse. Let  $X$ be not dicritical and let prove that
$f$ is not dicritical.
Let denote by $\tilde{f}$ and $\tilde{X}$ (respectively)
 the blow-up of the map $f$ and of the field $X.$
We have to
prove that if the exceptional divisor $D$ is invariant by $\tilde{X}$ then $\tilde{f}$ is tangential on $D.$
 We can suppose that:
$$
\tilde{X}=A(x,y)\frac{\partial}{\partial x}+B(x,y)\frac{\partial}{\partial y},
$$
and  $D=\{l(x,y)=x=0\}.$ The invariance of $D$  is equivalent to the fact that $x$ divides  $A(x,y).$ Let be
$$
T=\max\{s\in \mathbb{N}\mid x^s\mid A(x,y)\},
$$
i.e. $A(x,y)=x^T a(x,y).$ So we have:
$$
X=x^T a(x,y)\frac{\partial}{\partial x}+B(x,y)\frac{\partial}{\partial y},
$$
with $x \nmid B(x,y)$ i.e. $B(0,y)=y^k+ \cdots .$ Now let be $\bar{X}:= \pi^\ast (X)$ and observe that this field has the following structure:
$$
\bar{X}= x^\alpha a(x,y)\frac{\partial}{\partial x} + x^\beta B(x,y)\frac{\partial}{\partial y},
$$
where $a(x,y)$ and $B(x,y)$ are the previous ones and $\alpha > \beta.$
By Lemma \ref{lemma: marsed} we know that $\exp(\bar{X})= \tilde{f}$  and so we
can reconstruct the map by formula \eqref{eq:formal}.

We find $\bar{X}^j.x=x^\alpha(\cdots)$ for all $j.$
On the other hand when we compute $ \bar{X}^i.y,$ we find a structure of the
type $\bar{X}^i.y= x^\alpha (\cdots)+ x^{i\beta}(\cdots).$
Then the lowest power of $x$ appears in the term $\bar{X}.y$ and it is
$x^{\beta}y^k.$ So the order of $\tilde{f}$ on $D$ is $ \min(\alpha, \beta)=\beta$ and then
$$
\frac{l\circ \tilde{f} - l}{l^T}=\frac{\tilde{f}_1 - x}{x^\beta}
= x^{\alpha-\beta}(\cdots) \equiv 0 \hbox{\hspace{0.5cm}mod $\mathcal{I}(S)_p$.}
$$
\end{proof}
In this setting we have
\begin{propos}
Let $f\in \Phi_{\geq 2}(\mathbb{C}^2,0)$ be a dicritical holomorphic map tangent to the identity in $\C^2.$
Then there exist infinitely many RP curves.
\end{propos}
\begin{proof}
Since $f$ is dicritical then, by Proposition \ref{prop:dicritico}, $X$ is dicritical.
By Proposition \ref{prop:inf separ} there exist infinite separatrices and then by Proposition \ref{propos:sepfor}
we get the assertion.
\end{proof}

\end{document}